\documentclass[11pt]{article}
\usepackage{graphicx}
\usepackage{amssymb}

\title{On the Unboundedness of the First Eigenvalue of the Laplacian for $G$-Invariant Metrics}
\author{Paul Cernea}

\begin{document}
\maketitle

\begin{abstract}
In this note we partially answer a question posed by Colbois, Dryden, and El Soufi.  Consider the space of constant-volume Riemannian metrics on a connected manifold $M$ which are invariant under the action of a discrete Lie group $G$.  We show that the first eigenvalue of the Laplacian is not bounded above on this space, provided $M = S^n$, $G$ acts freely, and $S^n/G$ with the round metric admits a Killing vector field of constant length, or provided $M \neq T^n$ is a compact Lie group, and $G$ is a discrete subgroup.
\end{abstract}

Let $G$ be a Lie group acting on a smooth orientable closed connected manifold $M$ of dimension $n> 2$.  Given a Riemannian metric $g$ on $M$, let $\nabla$ be its Riemannian connection, and $\triangle=\nabla^i\nabla_i$ its Laplacian.  This Laplacian has eigenvalues $\lambda>0$ and associated eigenfunctions $u\in C^\infty(M)$ which are by definition solutions to

\begin{equation}
\triangle u +\lambda u = 0
\end{equation}

with $u$ not identically zero.  The set of eigenvalues is discrete, and its only limit point is $\infty$, so we can define $\lambda_1(g)$ as the lowest eigenvalue.  Define a scale-invariant functional,

\begin{equation}
g \longmapsto \Lambda_1(g) = \lambda_1(g)\textnormal{vol}^{2/n}(M,g)
\end{equation}

In \textbf{[CDE]} the following question is posed (\textit{Remark 3.1}):

\medskip

For $G$ a discrete group, can we construct $G$-invariant metrics such that the functional $\Lambda_1 = \lambda_1\textnormal{vol}^{2/n}$ becomes arbitrarily large?  Note that we are \textit{not} requiring the eigenfunctions to be $G$-invariant.

\medskip

Suppose $G$ is just the trivial group.  If $n=2$, Hersch showed that $\Lambda_1(g) \le 8\pi$ on $M = S^2$.  More generally, P. Yang and S. T. Yau showed that if $M$ has genus $\gamma$, then $\Lambda_1(g) \le 8\pi(1+\gamma)$.  (See [\textbf{SY}] for a discussion of these results.)  On the other hand, if $n >2$, it is well-known that $\Lambda_1$ can be arbitrarily large for any choice of $M$ (see [\textbf{CD}]).

\medskip

In this paper, we will consider nontrivial $G$ actions.  We will see that $\Lambda_1$ is unbounded for $M = S^n$ with certain values of $n$ and certain $G$-actions.  We will also see that $\Lambda_1$ is unbounded for $M \neq T^n$ a compact Lie group, and $G$ a discrete subgroup $G \leq M$ acting on the left by isometries, as a corollary to a theorem of Urakawa. 

\section{Spheres.}

Let

\begin{equation}
\bar g_{ij} = g_{ij} + tY_iY_j
\end{equation}

\textbf{Lemma.}  Suppose that $Y^i$ is a Killing vector field (i.e. $\nabla_iY_j+\nabla_jY_i=0$) of constant length $(|Y|\equiv \textit{const.})$.  Then the Laplacian $\bar \triangle$ of $\bar g$ satisfies

\begin{equation}
\bar \triangle u = \triangle u-\frac{t}{1+t|Y|^2}Y^iY^j\nabla_i\nabla_ju
\end{equation}

\textbf{Proof.}  First,

\begin{equation}
\bar g^{ij} = g^{ij} - \frac{tY^iY^j}{1+t|Y|^2}
\end{equation}

Next,

\begin{equation}
\bar \Gamma^k_{ij} = \frac{1}{2}\left(g^{kl}-\frac{tY^kY^l}{1+t|Y|^2}\right)\left(\left[g_{il,j}+g_{jl,i}-g_{ij,l}\right]+t\left[\left(Y_iY_l\right)_{,j}+\left(Y_jY_l\right)_{,i}-\left(Y_iY_j\right)_{,l}\right]\right)
\end{equation}

\begin{equation}
= \frac{1}{2}\left(g^{kl}-\frac{tY^kY^l}{1+t|Y|^2}\right)\left(g_{ls}\Gamma^s_{ij}+t\left[\left(Y_iY_l\right)_{,j}+\left(Y_jY_l\right)_{,i}-\left(Y_iY_j\right)_{,l}\right]\right)
\end{equation}

Let $K_{ij}=\nabla_iY_j+\nabla_jY_i$ and $C_{ij}=\nabla_jY_i-\nabla_iY_j=Y_{i,j}-Y_{j,i}$.  Then

\begin{equation}
\bar \Gamma^k_{ij}= \frac{1}{2}\left(g^{kl}-\frac{tY^kY^l}{1+t|Y|^2}\right)\left(g_{ls}\Gamma^s_{ij}+t\left[K_{ij}Y_l+2\Gamma^s_{ij}Y_sY_l+Y_iC_{lj}+Y_jC_{li}\right]\right)
\end{equation}

Thus

\begin{equation}
\bar \Gamma^k_{ij}=\Gamma^k_{ij}-\frac{tY^kY_s}{1+t|Y|^2}\Gamma^s_{ij}+\frac{t}{2}K_{ij}Y^k+t\Gamma^s_{ij}Y_sY^k+\frac{t}{2}g^{kl}\left(Y_iC_{lj}+Y_jC_{li}\right)-\frac{t^2K_{ij}Y^k|Y|^2}{2(1+t|Y|^2)}
\end{equation}

\begin{equation}
-\frac{t^2\Gamma^s_{ij}Y_sY^k|Y|^2}{1+t|Y|^2}-\frac{t^2Y^kY^l\left(Y_iC_{lj}+Y_jC_{li}\right)}{2(1+t|Y|^2)}
\end{equation}

\medskip

Now suppose $Y$ is a Killing vector field.  Then $K_{ij}=0$, and taking its trace, div$(Y)=0$.  Moreover, $C_{ij}=2\nabla_jY_i$, and

\begin{equation}
Y^i\nabla_iY^j = -Y^i\nabla^jY_i = -\frac{1}{2}\nabla^j|Y|^2=0
\end{equation}  

Hence $Y^lY_iC_{lj}=0$ and

\begin{equation}
\bar \Gamma^k_{ij}=\Gamma^k_{ij}-\frac{tY^kY_s}{1+t|Y|^2}\Gamma^s_{ij}+t\Gamma^s_{ij}Y_sY^k+\frac{t}{2}g^{kl}\left(Y_iC_{lj}+Y_jC_{li}\right)-\frac{t^2\Gamma^s_{ij}Y_sY^k|Y|^2}{1+t|Y|^2}
\end{equation}

However,

\begin{equation}
-\frac{tY^kY_s}{1+t|Y|^2}\Gamma^s_{ij}+t\Gamma^s_{ij}Y_sY^k-\frac{t^2\Gamma^s_{ij}Y_sY^k|Y|^2}{1+t|Y|^2}=0
\end{equation}

so

\begin{equation}
\bar \Gamma^k_{ij}=\Gamma^k_{ij}+\frac{t}{2}g^{kl}\left(Y_iC_{lj}+Y_jC_{li}\right)
\end{equation}

Let 

\begin{equation}
A^k_{ij}=\bar\Gamma^k_{ij}-\Gamma^k_{ij}=tg^{kl}\left(Y_i\nabla_jY_l+Y_j\nabla_iY_l\right)
\end{equation}

Then

\begin{equation}
\bar \triangle u = \triangle u-\frac{t}{1+t|Y|^2}Y^iY^j\nabla_i\nabla_ju -A^k_{ij}g^{ij}\nabla_ku+\frac{t}{1+t|Y|^2}A^k_{ij}Y^iY^j\nabla_ku
\end{equation}

But since

\begin{equation}
A^k_{ij}g^{ij}=A^k_{ij}Y^iY^j=0
\end{equation}

we obtain

\begin{equation}
\bar \triangle u = \triangle u-\frac{t}{1+t|Y|^2}Y^iY^j\nabla_i\nabla_ju
\end{equation}

as desired.  $\blacksquare$

Now rescale $Y$ so that $|Y|=1$. Then we have

\begin{equation}
\bar \triangle u = \triangle u-\frac{t}{1+t}Y^iY^j\nabla_i\nabla_ju
\end{equation}

Also, the volume transforms as $\bar\textnormal{vol}^{2/n}(M,g) = (1+t)^{1/n}\textnormal{vol}^{2/n}(M,g)$ since

\begin{equation}
\frac{\partial}{\partial t}d\mu = \frac{1}{2}\bar g^{ij}Y_iY_j d\mu=\frac{1}{2(1+t)}d\mu
\end{equation}

\medskip

Now specialize to the case where $g$ is the round metric on $M=S^n$ and $Y$ is a unit-length Killing field on $(M,g)$ with the property that $Y$ is the lift of a Killing field on $M/G$.  Then $\bar g$ is a one-parameter family of $G$-invariant metrics.

\medskip

Let $u$ be a first eigenfunction of $\triangle$ on $S^n$.  Any such $u$ is just a coordinate projection $x^1,x^2,x^3,\ldots$ or $x^{n+1}$ restricted to $S^n$ (or a linear combination thereof), and is characterized by

\begin{equation}
\nabla_i\nabla_ju = -g_{ij}u
\end{equation} 

Then

\begin{equation}
\bar \triangle u = \left(-n +\frac{t}{1+t}\right)u
\end{equation}

It follows that  $n-t/(1+t)$ is an eigenvalue of $\bar \triangle_t$.  But that's not all.  Let $\lambda_k(t)$ be the $k$th smallest eigenvalue of $\bar g(t)$, counted with multiplicity.  Since $\lambda_k(t)$ is continuous in $t$ (see [\textbf{BU}]), it follows that $n-t/(1+t)$ is the \textit{first} eigenvalue of $\bar\triangle_t$.

\medskip

To see this, let $S$ be the set of $t\in (-1,\infty)$ such that $\lambda_1(t)=n-t/(1+t)$.  We will show that $S$ is open and closed in $(-1,\infty)$, as well as nonempty.  Obviously, $0\in S$.  Next, by continuity of $\lambda_k(t)$, $S$ is closed.  Finally, let $\tau \in S$, and suppose that $\lambda_1(t) < n-t/(1+t)$ for $t$ close to $\tau$.  Then there exists a maximal natural number $m>n+1$ such that $n-t/(1+t) =\lambda_m(t)$.  But $\lambda_m(\tau) > n-\tau/(1+\tau)$ which contradicts the continuity of $\lambda_m(t)$ at $\tau$.  So $n-t/(1+t)$ is indeed the first eigenvalue of $\bar\triangle_t$ for all $t\in (-1,\infty)$.

\medskip

  Thus

\begin{equation}
\Lambda_1(t) = \left(n-\frac{t}{1+t}\right)(1+t)^{1/n}\textnormal{vol}^{2/n}(M,g(0))
\end{equation}

Since $n>1$, we see that $\Lambda_1 \rightarrow\infty$ as $t\rightarrow \infty$.  This proves unboundedness.

\medskip

Now it only remains to find examples of $S^n$ and $G$ for which $S^n/G$ admits a nonvanishing Killing field.  The Euler characteristic obstructs the existence of a nonvanishing vector field.  Even-dimensional spheres, then, are at once ruled out, since the Euler characteristic $\chi(S^{2m})$ is nonzero.  The Euler characteristic of odd closed manifolds, however, vanishes by Poincar\'e Duality.

\medskip

Let us make the following definition.  A Riemannian manifold $(M,g)$ is called a \textit{Sasakian manifold} if $M$ has a unit length Killing vector field $Y$ such that, for any vector fields $A$ and $B$ on $M$ we have

\begin{equation}
R(A,Y)B = g(Y,B)A - g(A,B)Y
\end{equation}

It is a classical result of Sasaki that all $3$-dimensional spherical space forms are Sasaki.  Thus examples of $G$ include any cyclic group $\mathbb{Z}_m$ (whose actions result in the Lens spaces).  More generally, $G$ can be any finite subgroup of $SO(4)$ acting freely by rotations on $S^3$.  These $G$ give rise to the so-called \textit{spherical 3-manifolds}.  By Grigori Perelman's proof of the \textit{elliptization conjecture}, these are all the possible fundamental groups of discrete quotients of $S^3$.  Such a $G$ is either cyclic, or a central extension of a dihedral, tetrahedral, octahedral, or icosahedral group by a cyclic group of even order.  For example, if $G$ is

\begin{equation}
\left<a,b | (ab)^2 = a^3 = b^5\right>
\end{equation}

then $S^3/G$ is the Poincar\'e homology sphere.

\medskip

Further higher-dimensional examples include $G = \mathbb{Z}_2$ acting on $S^{4m-1}$ since

\begin{equation}
\mathbb{R}P^{4m-1} = \frac{Sp(m)}{Sp(m-1)\times \mathbb{Z}_2}
\end{equation}

is Sasakian.  For these and more results on Sasakian manifolds, refer to [\textbf{BGM}].

\section{Compact Lie Groups.}

In [\textbf{U}] we have the following (\textit{Theorem 4}):

\medskip

\textit{Let $M$ be a compact connected Lie group.  We assume $M$ has nontrivial commutator subgroup.  That is, the commutator Lie subalgebra $\mathfrak{m}_1$ of $\mathfrak{m}$ is not zero.  Then there exists a family of left-invariant Riemannian metrics $g(t)$ ($0<t<\infty$) on $M$ such that}

\begin{equation}
\lim_{t\rightarrow\infty}\lambda_1(g(t))=\infty
\end{equation}

\begin{equation}
\lim_{t\rightarrow 0}\lambda_1(g(t))=0
\end{equation}

\textit{and $\textnormal{vol}(M,g(t))$ is constant in $t$.}

\medskip

To say that $\mathfrak{m}_1=0$ is to say that $\mathfrak{m}$ is abelian, and so $M$ is abelian, i.e. $M$ is a torus $T^n$, being compact.  Therefore, we have

\medskip

\textbf{Corollary.}  If $M \neq T^n$ is a compact connected Lie group, and $G$ is a discrete subgroup, then $\Lambda_1$ is unbounded among the $G$-invariant metrics on $M$, where $G$ acts by left-multiplication.

\medskip

\textbf{Proof.}  Since $M \neq T^n$, the existence of the family $g(t)$ from Theorem 4 is guaranteed to us.  Furthermore, all those metrics are $G$-invariant, since they are left-invariant.  $\blacksquare$

\section{References.}

\textbf{[CD]}  Bruno Colbois and J. Dozdiuk, \textit{Riemannian Metrics with Large $\lambda_1$}, Proc. Amer. Math. Soc., 122(3):905-906, 1994.

\textbf{[BGM]}  Charles P. Boyer, Krzysztof Galicki, and Benjamin M. Mann, \textit{Quaternionic Geometry and 3-Sasakian Manifolds}, Proceedings of the Meeting on Quaternionic Structures in Mathematics and Physics, Trieste (1994)

\textbf{[BU]}  Shigetoshi Bando and Hajime Urakawa, \textit{Generic Properties of the Eigenvalue of the Laplacian For Compact Riemannian Manifolds}, Tohoku Math. Journ. 35, 155-172 (1983)

\medskip

\textbf{[CDE]}  Bruno Colbois, Emily Dryden, and Ahmad El Soufi, \textit{Extremal G-Invariant Eigenvalues of the Laplacian of G-Invariant Metrics}, available at arXiv:math/0702547    

\medskip

\textbf{[SY]}  R. Schoen and S. T. Yau, \textit{Lectures on Differential Geometry, Volume 1}, pp. 135-136, (c) 1994. 

\medskip

\textbf{[U]}  Hajime Urakawa, \textit{On the Least Positive Eigenvalue of the Laplacian for Compact Group Manifolds}, J. Math. Soc. Japan, Vol 31, No. 1, 1979.

\end{document}